\documentclass[letterpaper,12pt]{article}
\usepackage{amsmath,amssymb,amsthm}
\usepackage{graphicx}
\usepackage{hyperref}
\usepackage{subfig}

\newtheorem{theorem}{Theorem}

\newtheorem{proposition}[theorem]{Proposition}
\newtheorem{lemma}[theorem]{Lemma}
\newtheorem{corollary}[theorem]{Corollary}

\theoremstyle{remark}

\newtheorem*{remarks}{{\bf Remarks}}

\newcommand{\I}[1]{I_{#1}}
\newcommand{\s}[1]{[ #1]}

\newcommand{\C}{{\bf C}}
\newcommand{\Q}{{\bf Q}}
\newcommand{\R}{{\bf R}}
\newcommand{\T}{{\bf T}}
\newcommand{\Z}{{\bf Z}}

\newcommand{\RU}{R(U)}
\newcommand{\RUn}{R(U_n)}
\newcommand{\im}{\operatorname{Im}}

\begin{document}
\title{Origami rings}
\author{Joe Buhler \thanks{Center for Communications Research, La Jolla, CA 92121 ({\tt buhler@ccrwest.org})} \and Steve Butler\thanks{Department of Mathematics, UCLA,
Los Angeles, CA 90095 ({\tt butler@math.ucla.edu}).\newline  This work was done with support of an NSF Mathematical Sciences Postdoctoral Fellowship.} \and Warwick de Launey\thanks{Center for Communications Research, La Jolla, CA 92121 ({\tt warwickdelauney2010@earthlink.net})} \and Ron Graham \thanks{Department of Computer Science and Engineering, University of California, San Diego,
La Jolla, CA 92093 ({\tt graham@ucsd.edu}).}}
\date{\empty}
\maketitle

\begin{abstract}
Motivated by a question in origami, we consider sets of
points in the complex plane constructed in the following
way.  Let $L_\alpha(p)$ be the
line in the complex plane through~$p$ with angle~$\alpha$ (with respect to the real axis).  Given a fixed collection $U$ 
of angles, let $\RU$ be the points that can be obtained by
starting with $0$ and~$1$, and then recursively adding intersection 
points of the form $L_\alpha(p) \cap L_\beta(q)$, where $p, q$ 
have been constructed already, and $\alpha, \beta$ are distinct
angles in~$U$.

Our main result is that if $U$ is a group with at least three
elements, then $\RU$ is a subring of the complex plane, i.e.,
it is closed under complex addition and multiplication.
This enables us to answer a specific question about origami
folds: if $n \ge 3$ and the allowable angles are the $n$
equally spaced angles $k\pi/n$, $0 \le k < n$, then
$\RU$ is the ring $\Z[\zeta_n]$ if $n$ is prime,
and the ring $\Z[1/n,\zeta_{n}]$ if $n$ is not prime,
where $\zeta_n := \exp(2\pi i/n)$ is a primitive $n$-th root
of unity.
\end{abstract}

\bigskip


\bigskip


\bigskip

\section{Introduction}
Origami constructions start with a flat sheet of paper and, through
a series of folds, arrive at a three-dimensional figure.  Some of
the folds produce pleats, whereas others produce needed 
reference points.  More information about the techniques of 
origami design with some amazing origami figures can be found in Lang \cite{lang},
and some mathematical aspects of origami can be found in \cite{alperin} 
and~\cite{alperinlang}.

The sets of reference points that can be constructed under various 
assumptions has been studied for theoretical reasons --- e.g., 
solving quartic equations \cite[pp. 285--291]{demaine} --- or for more
practical reasons --- e.g., approximating desired reference points
\cite{langR}.

\begin{figure}[hfbt]
\centering
\includegraphics[scale=1.15]{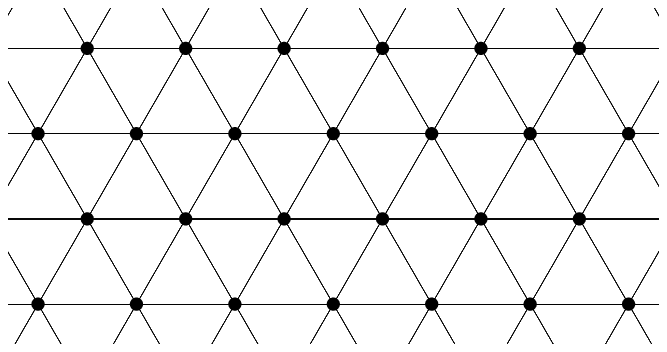}
\caption{The hexagonal lattice $\Z[\zeta_3]$. 
}
\label{fig:n3}
\end{figure}

We will consider an idealized form of paper folding and determine
the set of reference points that can be constructed by folds in a fixed set of directions (i.e., angles), where new points are found as intersections of folds through points that have already been constructed.
As a basic example, the set of all points constructed from the points $0$ and $1$ with
angles chosen from $\{0,\pi/3,2\pi/3\}$ is the hexagonal lattice
$\Z[\zeta_3]$ in the complex plane~$\C$ (see Figure~\ref{fig:n3}), as the reader may enjoy
verifying.

For $u\neq 0$, the line through $0$ and $u$ determines an angle with the real axis.  
Since 
nonzero real multiples of $u$ determine the same angle as $u$, without loss of generality we may assume $u$ has absolute value~$1$, i.e., that $u$ is on the unit
circle $\T$ in the complex plane.  In addition, since $u$ and $-u$ 
also determine the same angle, we choose to think of angles as
elements of the quotient group $\T/\{\pm 1\}$ of the unit circle
by its subgroup $\{\pm 1\}$ of order~$2$.  
We will abuse notation in the usual 
way and also refer to elements of~$\T$ as angles.
Also, note that the $1\in \T/\{\pm 1\}$ corresponds to a horizontal line 
(parallel to the real axis).  
The reader could alternatively choose to identify this group with the group 
$\R/\pi\Z$ of real numbers in $[0,\pi)$ under addition modulo~$\pi$, 
where $\theta$ corresponds to $e^{i\theta}$ in the formulas to follow.
In this setting, the horizontal line has angle~$0$.

To simplify and symmetrize our constructions, we will always assume 
our collection $U$ of angles is a group, i.e., that it is a subgroup
of $\T/{\pm 1}$.  If $U$ is finite, then it is cyclic; let
$U_n$ denote the group of angles that is cyclic of order~$n$
generated by $e^{i\pi/n}$.  Thus the collection $\{0,\pi/3,2\pi/3\}$ 
mentioned above corresponds to~$U_3$.

If $p$ is a point in the plane and $u$ is an angle then
\[
L_u(p) := \{ p+ru: r \in \R\}
\]
is the line through $p$ with angle~$u$.  We define the intersection
\[
\I{u,v}(p,q) := L_u(p) \cap L_v(q)
\]
to be the unique point on the intersection of the lines $L_u(p)$ and $L_v(q)$,
when $u$ and~$v$ are distinct angles.  Note that the slopes
of the two lines are distinct, so the intersection exists and is unique.

In general, given an initial (or seed) set $S$ lying in a set $X$ on which a collection 
of binary operators acts, it is natural to recursively define $S_0 = S$,
and $S_{n+1}$ to be the set obtained by applying all operators to
(pairs of) elements of $S_n$.  The union of the $S_n$ is the set
generated by the seed points $S$ by iterating the action of the operators,
this can also be defined as the smallest subset of~$X$ that contains~$S$
and is closed under the given set of operators.

We are interested in the case where the initial set $S$ consists of the two
points $0$ and $1$, and the operators are the set of $I_{u,v}$ where $u$ and~$v$ range
over distinct angles in~$U$.  
We let
$\RU$ denote the smallest subset of the complex plans that contains $\{0,1\}$
and is closed in the sense that if $p$ and $q$ are in the set, and $u$ and
$v$ are distinct elements of~$U$, then $\I{u,v}(p,q)$ is also in the set.

Our main result is the (perhaps slightly surprising) fact that $\RU$ is
a subring of $\C$ if $U$ has at least 3 elements.  Moreover, the elements of $\RU$ can be
described explicitly.  

\begin{theorem}\label{thm:ring}
With the above notation, if $U$ has at least 3 elements then $\RU$ is a subring 
of $\C$,  consisting of all integral linear combinations of 
``monomials,'' which are arbitrary finite products of complex numbers of the form
\[
\frac{1-u^2}{1-v^2}
\]
where $u$ and~$v$ are in $U$.  (Note that $u^2 \in \T$ is independent of the 
choice of representative of the class of $u \in \T/\{\pm 1\}$.)
\end{theorem}

This enables us to give a precise description of the ring when $U$ is finite.

\begin{theorem}\label{thm:finite}
Let $n \ge 3$.  If $n$ is prime, then $\RUn = \Z[\zeta_{n}]$ is 
the cyclotomic integer ring.  If $n$ is not a prime then
$\RUn= \Z[1/n,\zeta_{n}]$ is the cyclotomic integer ring localized at (the
primes dividing)~$n$.
\end{theorem}

This generalizes the $n = 3$ result mentioned above, but, as we will
explain later, that case is an outlier: for $n > 3$, $\RUn$ is a dense
subset of the complex plane.

\section{Plane intersections}
\label{sec:intersections}
The goal of this section is to describe basic facts about the intersection
point $\I{u,v}(p,q)$, both from algebraic and geometric perspectives.

Let $u,v$ be distinct angles, and consider the lines $L_u(p)$ and $L_v(q)$.  The intersection of these lines is the point $\I{u,v}(p,q)$ that is the unique solution to 
\[
z = p+ru = q+sv,
\]
for $r,s$ real numbers.  Since $s$ is real then 
the imaginary part of $s = v^{-1}(p-q+ru)$ must equal zero.  This gives an
equation that can be solved for~$r$, yielding
\[
r = \frac{\im((q-p)/v)}{\im(u/v)}.
\]
It is convenient to introduce the notation
\[
\s{x,y} = x y^* - x^* y = 2i \: |y|^2 \:  \im(x/y) 
\]
where $x^*$ denotes the complex conjugate of $x$, and $|y|^2 = y^*y$.
Note that $\s{\cdot,\cdot}$ is additive,
real-linear, and antisymmetric.  We can now rewrite the equation for~$r$ in the form
\[
r = \frac{\s{q-p,v}}{\s{u,v}}.
\]
This gives
\[
\I{u,v}(p,q) = p + \frac{\s{q-p,v}}{\s{u,v}} \: u =  
\frac{\s{u,v} p + \s{q,v} u - \s{p,v} u} {\s{u,v}}.
\]
Substituting the definition and doing some algebraic juggling leads
to a fundamental formula:
\begin{equation}
\label{eqn:Iform}
\I{u,v}(p,q) = 
\frac{up^*v-u^*pv-vq^*u+v^*qu} {\s{u,v}}
= \frac{\s{u,p}}{\s{u,v}} v + \frac{\s{v,q}}{\s{v,u}} u.
\end{equation}
A number of basic facts follow from this formula.

\begin{proposition}\label{prp:basicIfacts}
Let $p,q$ be points in the plane, and $u,v$ be pairwise distinct angles.
\begin{itemize}
\item (Symmetry)   $\I{u,v}(p,q) = \I{v,u}(q,p)$
\item (Reduction)  $\I{u,v}(p,q) = \I{u,v}(p,0) + \I{v,u}(q,0)$
\item (Projection)  $\I{u,v}(p,0)$ is a projection of~$p$ onto the line 
$\{rv : r \in \R\}$ in the direction~$u$
\item (Linearity)  $\I{u,v}(p+q,0) = \I{u,v}(p,0)+\I{u,v}(q,0)$ and, for real $r$,
$\I{u,v}(rp,0) = r \: \I{u,v}(p,0)$
\item (Convexity)  $\I{u,v}(p,q)$ has the form $Ap+Bq$ where $A$ and $B$ are
real-linear maps of the complex plane that satisfy 
$A+B = 1_\C$, where $1_\C$ is the identity map on the plane
\item (Rotation)  For $w \in \T$, $w\I{u,v}(p,q) = \I{wu,wu}(wp,wq)$.
\end{itemize}
\end{proposition}

\begin{remarks}
\begin{enumerate}
\item The projection of $p$, onto a line $L$ in the direction $u$, is the
point of the form $p+ru$ that lies on the line~$L$.
\item  Real-linear maps of the plane can be written in two useful forms:
as maps taking $p$ to $p' = ap+bp^*$ for complex constants $a,b$, or as
maps taking $p = x+iy$ to $p' = x'+iy'$ written in terms of two by two
matrices:
\[
\begin{pmatrix}
x' \\ 
y'
\end{pmatrix}
= 
\begin{pmatrix}
a & b \\
c & d 
\end{pmatrix}
\begin{pmatrix}
x \\ 
y
\end{pmatrix},
\]
where $a,b,c,d$ are real.  
\item The use of the term ``Convexity'' is merely by vague analogy with the 
operation $tx+(1-t)y$ taking elements $x,y$ of a real vector space to 
a convex combination on the line between them, where $t \in [0,1]$.
\end{enumerate}
\end{remarks}

The verifications of the facts in the proposition from the identity
(\ref{eqn:Iform}) are straightforward algebraic exercises which we leave
to the reader.  However, these facts
can also be verified geometrically, as we now briefly sketch.

\begin{figure}[hftb]
\centering
\hfil 
\subfloat[Law of parallelograms]{\includegraphics[scale=1.7]{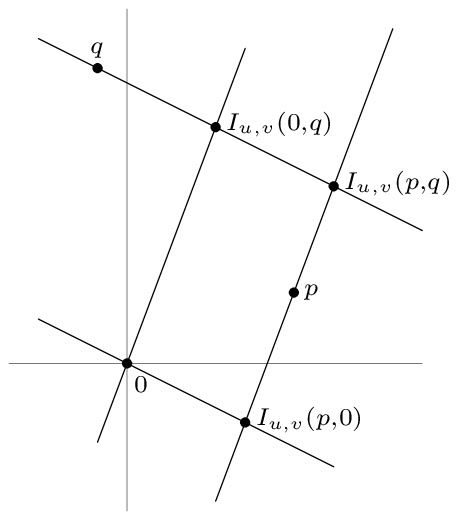}\label{fig:para}}
\hfil\hfil 
\subfloat[Projection]{\includegraphics[scale=1.7]{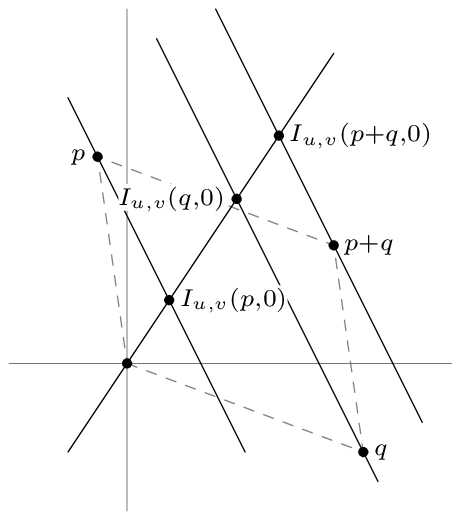}\label{fig:proj}}
\hfil\hfil
\caption{Properties of $\I{u,v}(p,q)$}
\end{figure}

Symmetry follows from the fact that both sides of the 
equation solve the same intersection problem.  Reduction is the Law of Parallelograms 
(see Figure~\ref{fig:para}) and Symmetry.  Projection follows from another straightforward
diagram.
Linearity follows from Projection and the fact that projections are 
linear (see Figure~\ref{fig:proj}).  Convexity follows from Reduction, 
the fact that each term in the Reduction formula is a projection, and the
fact that if $p = q$ then $\I{u,v}(p,q) = p$.  (Note that the linear maps
$A,B$ have rank 1, since they are projections.)  Rotation expresses the
fact that if the points and angles of an intersection problem are all 
rotated by an angle~$w$ then the result is rotated by~$w$.

\section{Closure under addition and multiplication}
\label{sec:thm1}
Recall for a fixed group $U$ of angles, $\RU$ denotes
the smallest set of complex numbers that contains $0$ and $1$ and has the
property that if $u,v$ are distinct elements of~$U$, and $p,q$ are elements
of $\RU$, then $\I{u,v}(p,q)$ is in~$\RU$.  Theorem~\ref{thm:ring} asserts
that if $U$ has size at least three then $\RU$ is a subring of~$\C$, i.e., that $\RU$ is closed under addition
and multiplication.

This will be proved in a series of lemmas.

\begin{lemma}\label{lem:addclos}
$\RU$ is closed under addition.
\end{lemma}

\begin{figure}[hftb]
\centering
\includegraphics[scale=1.25]{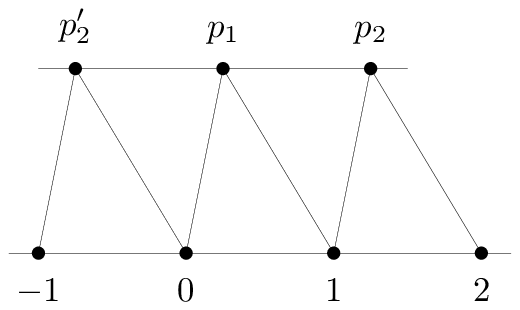}
\caption{$-1$ and $2$.}
\label{fig:pts}
\end{figure}

\begin{proof}  First, note that $2$ is in~$\RU$.  Indeed,
if $u,v$ are distinct nonzero angles in~$U$ then construct
the following points:
\[
p_1 = \I{u,v}(0,1), \quad  p_2 = \I{u,1}(1,p_1), \quad p_3 = \I{1,v}(0,p_2).
\]
Some slightly involved algebra shows that $p_3 = 2$, but this
is also geometrically obvious as a quick glance at Figure~\ref{fig:pts} shows.

Moreover, if $p$ is in $\RU$ then so is $p+1$.   Indeed, apply the
steps that were used to construct~$p$ starting from 0 and~1, but instead
starting at the points 1 and~2; the result is~$p+1$.  Similarly,
if $p$ and~$q$ are in $\RU$ then so is their sum: repeating the
steps used to construct $q$, starting from $p$ and $p+1$, gives
$p+q$.
\end{proof}

\begin{lemma}\label{lem:invclos}
$\RU$ is closed under taking additive inverses.
\end{lemma}

\begin{proof}
First, note that $-1$ is in~$\RU$.  Indeed, if $p_1 = \I{u,v}(0,1)$ as above,
and $p_2' = \I{v,1}(0,p_1)$ then it is easy to verify geometrically or algebraically 
that $-1 = \I{1,u}(0,p_2')$ (see Figure~\ref{fig:pts}).
If $p$ is in $\RU$ then $-p$
can be constructed by applying the steps used to construct~$p$, except
that we start at $0$ and~$-1$.
\end{proof}

A point $p$ in $\RU$ is said to be a \emph{monomial} (or, if the emphasis 
is needed, a $U$-monomial) if it can be produced starting at~$1$
using only intersections of the form $\I{u,v}(p,0)$, i.e., intersections
in which the second line is required to go through the origin.  Thus
a monomial lies in a sequence of the form
\[
p_1 = \I{u_1,v_1}(1,0), \quad p_2 = \I{u_2,v_2}(p_1,0), \quad p_3 = \I{u_3,v_3}(p_2,0), \quad \ldots
\]
If $p = p_n = \I{u_n,v_n}(p_{n-1},0)$ occurs at the $n$-th step, then $p$
is said to be a monomial of length at most~$n$.
In terms of the linear operators $A$ in the Convexity property,
this could be written
\[
p = A_n A_{n-1} \cdots A_2 A_1 1
\]
where $A_i$ is the linear operator occurring in the $i$-th step.

A monomial $p = I_{u,v}(1,0)$ of length 1 is said to be an \emph{elementary}
monomial.

\begin{lemma}\label{lem:mon}
$\RU$ is the set of (finite) integer linear combinations of monomials.
\end{lemma}

\begin{proof}
Since $\RU$ is closed under addition and contains~$-1$ it is obvious that $\RU$ contains
all linear combinations of monomials with coefficients in~$\Z$.

To prove that $\RU$ consists exactly of those elements, use induction
on the length of the constructions.   Specifically, we prove that elements
of $\RU$ constructed with at most $n$ intersection operations are linear
combinations of monomials of length at most~$n$.  The base cases are
obvious.  If $p$ and $q$ can be constructed with at most~$n$ operations,
then consider a point $z = \I{u,v}(p,q)$ that can be constructed in $n+1$
operations.  By the induction assumption, $p$ and~$q$ are sums of monomials
of length at most~$n$.
By Reduction and Linearity $z$ is a linear combination of terms of the form
$\I{u,v}(m,0)$ where $m$ is a monomial of length at most~$n$.  It follows
that $z$ is a linear combination of monomials of length at most~$n+1$.
\end{proof}

\begin{lemma}\label{lem:mul}
The product of two monomials is a monomial.
Any monomial is a product of elementary monomials.
The set $\RU$ is closed under multiplication.
\end{lemma}

\begin{proof}
For the first two statements, the key idea is the Rotation property
(or, that the set of operators is normalized by rotations).  The
identity
\[
\I{u,v}(1,0) \I{u',v'}(1,0) = r v \: \I{u',v'}(1,0) = 
r \: \I{vu',vv'}(v,0) = \I{vu',vv'}(rv,0) = 
\I{vu',vv'}(\I{u,v}(1,0),0)
\]
(where $r = \s{u,1}/\s{u,v} \in \R$)
shows that the product of two elementary monomials is a monomial of length~2,
and conversely.  Everything else is induction.

For instance, assume that the product of two monomials is a monomial when the 
sum of the lengths is at most~$n$.  Let $m'$ be a monomial of length at most~$a$, 
and $m$ be a monomial of length at most~$b$ where $a+b = n$.  Then consider
a product of monomials of length $a+1$ and~$b$:
\[
m \: \I{u,v}(m',0).
\]
Write $m = tw$ where $t$ is real and $w$ is on the unit circle.
Use Rotation and Linearity to get
\[
m \: \I{u,v}(m',0) = tw \: \I{u,v}(m',0) = 
\I{wu,wv}(twm',0) = \I{wu,wv}(mm',0).
\]
Since $mm'$ is a monomial by the induction assumption, it follows that the right-hand side
is a monomial.  A similar induction argument proves the converse.

The last statement of the lemma follows immediately from the first statement,
using the preceding lemma.
\end{proof}

This finishes the proof that $\RU$ is a subring of the complex numbers.

\begin{remarks}
\begin{enumerate}
\item The fact that our allowable angles $U$ formed a group was only
used in the last lemma, where we needed closure under multiplication.  
\item Many of the facts can be proved under weaker hypothesis, and these
results might be useful in other contexts.
For instance, if we consider more general operators $Ap+Bq$, as in
the Convexity property, applied recursively to a set of seed points,
then closure under addition is implied by (and therefore equivalent to) 
the fact that $2$ can be constructed.
In similar more general contexts, the property of being closed under
multiplication comes down to the set of operators being normalized by
rotations (multiplying by elements of~$\T$).
\end{enumerate}
\end{remarks}

\section{Ring generators}

Throughout this section let $U$ be a group of angles with 3 or more
elements, and let $U_n$ denote the cyclic group of order $n$ of angles 
generated by $e^{i\pi/n}$, for $n \ge 3$.  

From the proof of multiplicative closure, monomials have the form
\[
\frac{\s{u_1,1}}{\s{u_1,v_1}}
\frac{\s{u_2,u_1}}{\s{u_2,v_2}} \cdots
\frac{\s{u_n,u_{n-1}}}{\s{u_n,v_n}} \: v_n.
\]
If $u = e^{i\pi a}$, $v =  e^{i\pi b}$, and $w = e^{i\pi c}$ then
\[
\frac{\s{u,w}}{\s{u,v}} = \frac{\sin(\pi(a-c))}{\sin(\pi (a-b))}
\]
and we sometimes refer to these as ``sine-quotients.''  Thus
we can say that a monomial is a product of sine-quotients by an
element of~$U$.

The proof of multiplicative closure also showed that monomials are
products of elementary monomials.  
Note that an elementary monomial can be written, using $v = 1/v^*$, as
\[
\I{u,v}(1,0) = \frac{\s{u,1}}{\s{u,v}} \: v = \frac{u-u^*}{u(v^*)^2-u^*}
= \frac{1-u^2}{1-(u/v)^2}.
\]
As $u$ and $v$ range over all angles, this quotient ranges over all
quotients of the form $(1-u)/(1-v)$ where $u$ and $v$ are elements of
the group $V := \{u^2 : u \in U\} \subset \T$ obtained by taking 
representatives of elements of $U$, in the circle group, and squaring.  
(Note that $V$ is isomorphic to~$U$; indeed, the squaring map induces
an isomorphism $\T/\{\pm 1\} \simeq \T$, and $V$ is the image of~$U$
under this map.)

Combining all of these facts gives a proof of the basic assertion
made about $\RU$ in the introduction.

\begin{theorem}
Let $U$ be a group of angles with at least three elements, and let
$V \subset \T$ be the group of squares described above.  Then
$\RU$ is the subring of $\C$ consisting of linear combinations of
monomials, with coefficients in~$\Z$, where monomials are products of
elements of the form
\[
\frac{1-u}{1-v}
\]
for $u$ and~$v$ in~$V$, $v \ne 1$.
\end{theorem}

Recall that the cyclotomic field $\Q(\zeta_n)$ is the smallest subfield
of $\C$ that contains the primitive $n$-th root of unity $\zeta_n := e^{2\pi i/n}$.
It is a standard fact from algebra that $\Q(\zeta_n)$ is 
the set of rational functions in $\zeta_n$, over $\Q$, and that it
can be described even more simply as the set of
rational linear combinations of powers of~$\zeta_n$.
The group $U_n$ is generated by (the class of) $\zeta_{2n} = e^{i\pi/n}$,
and its group $V_n$ of squares is generated by $\zeta_n$.  The preceding
theorem immediately implies the following corollary.

\begin{corollary} Every element of $\RUn$ lies in the cyclotomic
field $\Q(\zeta_n)$.
\end{corollary}

\begin{remarks}
\begin{enumerate}
\item  The equality between the two formulas for $\I{u,v}(0,1)$
\[
\frac{\s{u,1}}{\s{u,v}} \: v = \frac{1-u^2}{1-(u/v)^2}
\]
is crucial in the above analysis, but it is also curious from a 
field-theoretic point of view.  The right hand
side is obviously in $\Q(\zeta_n)$.  However, the left-hand side is a product
of a power of $\zeta_{2n}$ and a sine-quotient lying in the totally real subfield of
$\Q(\zeta_{2n})$.
\item For later use, we remark that if $d$ is a divisor of~$n$ then
\[
\zeta_n^{n/d} = \zeta_d.
\]
\end{enumerate}
\end{remarks}

The sine-quotients that arise in the above algebra are familiar (at least
when $n$ is prime) from algebraic number theory as ``cyclotomic units''.
To analyze the integrality properties of these quotients we need to
recall some algebraic number theory.

The algebraic integers in the field $\Q(\zeta_n)$ form a ring
$\Z[\zeta_n]$ consisting of of integral linear combinations of powers
of $\zeta_n$.
The following facts can be found in many texts on algebraic number 
theory, such as \cite{cyclotomic}, and we include brief proofs for the convenience of the reader.

\begin{theorem}  Fix~$n > 3$, let $\zeta = \zeta_n$, and
assume throughout that $a$ and $b$ are nonzero modulo~$n$.
\begin{description}
\item[(a)] If $n$ is prime, then $(1-\zeta^a)/(1-\zeta^b)$ is an algebraic integer.
\item[(b)] If $n$ is non-prime, then $(1-\zeta^a)/(1-\zeta^b)$ has denominator
dividing~$n$.
\item[(c)] If $n$ is non-prime and $p$ is a prime divisor of~$n$,
then some product of an element of $\Z[\zeta_n]$ by a product
of quotients of the form $(1-\zeta^a)/(1-\zeta^b)$
is equal to $1/p$.
\end{description}
\end{theorem}

\begin{proof}
If $n$ is prime then $b$ is relatively prime to $n$, and we can
solve the equation $a+rn = bs$ in integers $r,s$.  Then 
\[
\frac{1-\zeta^a}{1-\zeta^b} = 
\frac{1-\zeta^{a+rn}}{1-\zeta^b} =
\frac{1-\zeta^{bs}}{1-\zeta^b}
\]
is obviously a polynomial in $\zeta$ and hence an algebraic integer.  This
proves~(a).

To prove the next two parts of the theorem we start with the identity
\[
X^n-1 = \prod_{k=0}^{n-1}(X-\zeta^k).
\]
Divide by $X-1$ and take the limit as $X \rightarrow 1$ to get
\begin{equation}
\label{eqn:enn}
n = \prod_{k=1}^{n-1}(1-\zeta^k).
\end{equation}
(We will refer to this as equation $(\ref{eqn:enn})_n$.)
It follows from $(\ref{eqn:enn})_n$ that $1-\zeta^k$ is a divisor of~$n$,
for all $k$ not divisible by~$n$, proving~(b).

To prove (c), note that it suffices to consider the cases $n = pq$
and $n = p^2$, where $p$ and~$q$ are primes.  Indeed, if $n$ is non-prime then
some divisor $d$ of $n$ has one of those forms, and we can confine ourselves to
looking at quotients $(1-\zeta^a)/(1-\zeta^b)$ where both $a$ and~$b$ are
multiples of $n/d$, which is equivalent to restricting to $n$ of those two forms.

If $n = pq$ then divide equation $(\ref{eqn:enn})_n$ by equations
$(\ref{eqn:enn})_p$ and $(\ref{eqn:enn})_q$ to get
$1 = \prod (1-\zeta^k)$ where $k$ ranges over integers relatively prime to~$n$.
It follows that $1-\zeta$ is a unit in $\Z[\zeta]$.  Dividing $(1-\zeta)^{p-1}$
by equation $(\ref{eqn:enn})_p$ gives
\[
\frac{u}{p} = \prod_{k=1}^{p-1} \frac{1-\zeta}{1-\zeta^{qk}}
\]
where $u$ is a unit, and we recall that $\zeta^q = \zeta_p$.  Thus $1/p$ is
a product of the claimed form.

If $n = p^2$ then careful use of equations $(\ref{eqn:enn})_{p^2}$ and 
$(\ref{eqn:enn})_p$ shows that
\[
\prod \frac{1-\zeta^{pk}}{1-\zeta^k} = \frac{1}{p^{p-1}}
\]
where the product is over $k$ up to $p^2$ and relatively prime to~$p$.
Multiplying by $p^{p-2}$ gives $1/p$ as desired.
\end{proof}

Our main theorem about $\RUn$ follows easily.

\begin{theorem}
If $n = p$ is prime, then $\RUn$ is equal to $\Z[\zeta_n]$.  
If $n$ is non-prime, then $\RUn$ is equal to $\Z[1/n,\zeta_n]$.
\end{theorem}

\begin{proof}
The ring $\RUn$ contains $-1$ and
\[
\frac{1-\zeta}{1-\zeta^{-1}} = -\zeta
\]
so it contains $\Z[\zeta_n]$.  If $n$ is prime then the first part of
the preceding theorem says that elementary monomials, and hence all
monomials, are in $\Z[\zeta_n]$, and it follows that
every element of $\RUn$ is contained in that ring
and we are done.  If $n$ is non-prime, then every element of $\RUn$ is
contained in $\Z[1/n,\zeta_n]$, by the second part of the preceding
theorem.  However, it is also easy to construct $1/n$ by using the 
third part of the preceding theorem, and hence any element of that ring, 
so that any element of $\Z[1/n,\zeta_n]$ is in $\RU$, and $\RUn = \Z[1/n,\zeta_n]$ 
as claimed.
\end{proof}

\begin{corollary} If $n > 3$ then $\RUn$ is dense in the complex plane.
\end{corollary}

\begin{proof}  
It is well-known in algebraic number theory that $\Z[\zeta_n]$ is dense
in the plane if $n \ge 4$, and the result follows easily.  However, it can
also be proved directly: denseness of the ring $\RUn$ is equivalent to 
showing that $0$ is a cluster point of $\RUn$, which in turn is equivalent
to finding a nonzero element of $\RUn$ inside the unit circle (since powers
of such an element converge to~$0$).   This is obvious for non-prime~$n$,
so we need only consider primes $n \ge 5$.  In this case
\[
\frac{1-\zeta}{1-\zeta^2} = \frac{1}{1+\zeta}
\]
is easily verified to be a unit lying inside the unit circle.
\end{proof}

\section{Concluding remarks}\label{sec:conclusion}
Motivated by a question about origami constructions for which the folds are limited
to a finite set of directions, we have been led to a general construction in
which  a group of angles~$U$ gives rise to a subring $\RU$ of the complex plane.  
When combined with standard results from the theory of cyclotomic fields,
the ring $\RU$, this gives an explicit description of $\RU$ in 
the motivating case in which~$U$ is finite.  

This may not be terribly useful in real-world
origami.  In the case of $U_n$, there is the problem about how to fold
the angles in the first place; in addition, although our proofs
are constructive they are likely to produce tedious series of 
folds --- in the real world, approximations with few folds would be 
much more valuable.  

However, the mathematical ideas that arise seem to lead to interesting
questions, and we mention several.

The set $R(U)$ can be a ring even if $U$ is only closed under multiplication
and not necessarily under inversion.  It is even a ring in the case
$U=\{0,\pi/4,\pi/2\}$.  Which sets of angles $U$ give rise to sets that
are closed under complex multiplication?

Which subrings of $\C$ are of the form $\RU$ for some group $U$ of angles?

Given a collection $\{A_i\}$ of $\R$-linear maps on the complex plane, define 
$B_i = 1_\C-A_i$, and consider the operators
\[
I_i(p,q) = A_ip+B_iq
\]
on the complex plane.  Given two seed points, recursively apply these
operators to build the subset of the plane of all points that are
constructable from the given set of operators.  (By Convexity, this 
generalizes~$\RU$.) What can be said about these sets?

\section*{Acknowledgments}
We thank Erik Demaine for originally suggesting this problem to us.

\end{document}